\documentclass[11pt]{amsart}
\usepackage{amssymb,amscd}

\begin{document}

\newtheorem{definition}{Definition}
\newtheorem{theorem}{Theorem}
\newtheorem{remark}{Remark}
\newtheorem*{lemma}{Lemma}
\newtheorem{corollary}{Corollary}
\newtheorem*{question}{Question}

\title{A short proof of  nonhomogeneity of the pseudo-circle}

\author{Krystyna~Kuperberg}
\address{Department of Mathematics and Statistics, Auburn University, Auburn, AL
36849, USA}
\email{kuperkm@auburn.edu}

\author{Kevin~Gammon}
\address{Department of Mathematics and Statistics, Auburn University, Auburn, AL
36849, USA}
\email{gammokb@auburn.edu}

\dedicatory{Dedicated to  James T.~Rogers,~Jr. on the occasion of his 65th birthday}

\begin{abstract} The pseudo-circle is known to be nonhomogeneous.  The original proofs of this fact were discovered independently by L.~Fearnley~\cite{Fearnleyhomogeneous} and J.T.~Rogers, Jr.~\cite{Rogers1}. The purpose of this paper is to provide an alternative, very short proof  based on  a result of D.~Bellamy and W.~Lewis~\cite{BellamyLewis}.
\end{abstract}

\maketitle

\section{Introduction}

A {\em pseudo-arc} is a hereditarily indecomposable, chainable continuum.   In 1948, E.E.~Moise~\cite{Moise} constructed a pseudo-arc as an indecomposable continuum  homeomorphic to each of its  subcontinua. Moise correctly conjectured that the  hereditarily indecomposable continuum given by B.~Knaster~\cite{Knaster}  in 1922 is a pseudo-arc. Also in 1948, R.H. Bing~\cite{Bing1} proved that  Moise's example is homogeneous.  In 1951, Bing~\cite{Bing2} proved that every hereditarily indecomposable chainable continuum is a pseudo-arc and that all pseudo-arcs are homeomorphic.  In 1959, Bing~\cite{Bing3} gave another characterization of the pseudo-arc: a homogeneous chainable continuum.

The history of many other aspects of the pseudo-arc can be found in  survey papers by W.~Lewis \cite{Lewis2} and \cite{Lewis1}.  

In 1951, Bing~ \cite{Bing2} described a  {\em pseudo-circle}, a planar  hereditarily indecomposable circularly chainable  continuum which separates the plane and whose every proper subcontinuum is a pseudo-arc.  It has been shown by  L.~Fearnley in~\cite{Fearnleyhomogeneous} and J. T.~Rogers, Jr. in ~\cite{Rogers1} that the pseudo-circle is not homogeneous.  Fearnley also proved that the pseudo-circle is unique~\cite{Fearnleyunique1} and \cite{Fearnleyunique2}. The fact that the 
pseudo-circle is not homogeneous also follows from more general theorems proved in~\cite{Hagopian}, \cite{KennedyRogers},   
\cite{Lewis3}, and \cite{Rogers2}.

This paper offers yet another, very short  proof, a consequence of  a result of D.~Bellamy and W.~Lewis~\cite{BellamyLewis}. Similarly as in~\cite{Rogers2}, an infinite covering space of a plane separating continuum is used.

\section{Preliminaries}\label{prelim}
Throughout the paper, a {\em continuum\/} will refer to a nondegenerate compact and connected metric space.  A continuum is {\em indecomposable} if it is not  the union of two proper subcontinua.  A  continuum is {\em hereditarily indecomposable} if every subcontinuum is also indecomposable.  
For a point $a$ in $X$, the {\em composant\/} $K(a)$ of $a $ in $X$ is the union of all proper subcontinua of $X$  containing $a$.  An indecomposable continuum contains uncountably many pairwise disjoint  composants, see~\cite{Kuratowski} Theorem 7,  page 212.

A topological space $X$ is {\em homogeneous} if for any two points in $X$ there is a homeomorphism of $X$ onto itself mapping one point onto the other.  

Let $C$ denote the pseudo-circle. We may assume that $C$ is contained in a planar annulus $A$ in such a way that the winding number of each circular chain in the sequence of crooked circular chains defining $C$ is one. Any homeomorphism $h:C \to C$  extends to a continuous map $f: A \to A$ of degree $\pm 1$. (First extend $h$ to a map  $\overline{U}\to A$ for some closed annular neighborhood $\overline{U}$ of $C$, then compose a retraction of $A$ onto $\overline{U}$ with this extension.)

Let  $\widetilde{A}$ be the universal covering space of $A$ with projection $p$. For any $\widetilde{x}\in \widetilde{A}$ and $\widetilde{y}\in  p^{-1}(f(p(\widetilde{x})))$ there is a map $\widetilde{f}$ such that the diagram 

$$
\begin{array}{ccc}
 &\widetilde{f} & \\
\widetilde{A} & \longrightarrow & \widetilde{A}\\
p\downarrow & & \downarrow p\\
A & \longrightarrow & A\\
& f &
\end{array}
$$

\noindent commutes and $\widetilde{f}(\widetilde{x})=\widetilde{y}$; see for example~\cite{Hu}, Theorem 16.3. Let $\widehat{A}$ be the disc that is a two-point compactification of  $\widetilde{A}$. Denote the two added points of the compactification  by $a$ and $b$.  The map $\widetilde{f}$ extends uniquely to a map $F: \widehat{A} \to \widehat{A}$.

Let $\widetilde{C}=p^{-1}(C)$, and let $P=\widetilde{C} \cup \{a,b\}$, 
a two-point compactification of $p^{-1}(C)$.  D.~Bellamy and W.~Lewis considered this set in~\cite{BellamyLewis}  and proved that $P$ is a pseudo-arc.

Denote  by $H$ the restriction of $F$ to $P$ and note that

\begin{enumerate}
\item either $H(a)=a$ and $H(b)=b$, or $H(a)=b$ and $H(b)=a$,
\item $\widetilde{f}(\widetilde{C})=\widetilde{C}$ and hence $H(P)=P$,
\item $\widetilde{f}_{|\widetilde{C}}$ is one-to-one.

\end{enumerate}

Thus 

\begin{lemma}  $H$ is a  homeomorphism from $P$ to $P$.\end{lemma}

\section{Proof of nonhomogeneity of the pseudo-circle}

\begin{theorem} The pseudo-circle is not homogeneous.
\end{theorem}

\begin{proof} Let $K(a)$ and $K(b)$ be the composants of $a$ and $b$, respectively, in the pseudo-arc $P$. Let $\widetilde{x}$ and $\widetilde{y}$ be two points in $P$ such that $\widetilde{x} \in (K(a)\cup K(b)) -\{a,b\}$ and $\widetilde{y}\in P-(K(a)\cup K(b))$. If $C$ were homogeneous, then there would be a 
homeomorphism $h:C\to C$ taking $x=p(\widetilde{x})$ onto $y=p(\widetilde{y})$. Then there would be a homeomorphism $H:P\to P$ as described in section~\ref{prelim} taking $\widetilde{x}$ onto $\widetilde{y}$.  This is not  possible since under every such homeomorphism, the set $K(a)\cup K(b)$ is invariant; the image of a composant is a composant.
\end{proof}

\noindent  {\bf Remark.}  It is not important for this proof that $K(a)$ and $K(b)$ are not the same set, but the authors are grateful to D.~Bellamy and W.~Lewis for showing that $K(a)$ and $K(b)$ were indeed different  composants.

\begin{theorem} If  for some $x$, the composant $K(a)$  intersects the fiber $p^{-1}(x)$, then it contains $p^{-1}(x)$.
\end{theorem}

\begin{proof} If $y \in p^{-1}(x)\cap K(a)$, then by the definition of a composant, there  is a proper subcontinuum $W$ of $P$ that contains both $a$ and $y$. Let $g:\widetilde{C}\to \widetilde{C}$ be a deck transformation such that $p^{-1}(x)=\{g^n(y)\}_{n\in \mathbb{Z}}$, $\mathbb{Z}$ being the set of integers. Denote by $G$ the extension of $g$ to $P$. The set $W_n=G^n(W)$ is a continuum containing $a$ and $g^n(y)$. Thus  $p^{-1}(x)\subset K(a)$.
\end{proof}

Note that Theorem 2 and Remark above imply that $p(K(a)-\{a\})\cap p(K(b)-\{b\})=\emptyset$.

\begin{question} Can the sets $p(K(a)-\{a\})$ and $ p(K(b)-\{b\})$ be used to classify  the composants of the pseudo-circle $C$?\end{question}

The authors would like to thank Jim Rogers, David Bellamy, and Wayne Lewis  for their comments.


\begin{thebibliography}{10}

\bibitem{Bing1} R.H. Bing, {\em A homogeneous indecomposable plane
    continuum\/}, Duke Math. J. 15 (1948), 729-742.
    
    \bibitem{Bing2} R.H. Bing, {\em Concerning hereditarily indecomposable
    continua\/}, Pacific J. Math. 1 (1951), 43-51.
    
   \bibitem{Bing3} R.H. Bing, {\em Each homogeneous nondegenerate chainable continuum is a pseudo-arc\/}, Proc. Amer. Math. Soc.  10  (1959),  345-346

\bibitem{BellamyLewis} D.P. Bellamy and W. Lewis, {\em An orientation reversing homeomorphism of the plane with invariant pseudo-arc\/}, Proc. Amer. Math. Soc. 114 (1992), 1145-1149.

\bibitem{Fearnleyunique1}  L. Fearnley, {\em  The pseudo-circle is unique\/},  Bull. Amer. Math. Soc.  75 (1969), 398-401.

\bibitem{Fearnleyhomogeneous} L. Fearnley, {\em The pseudo-circle is not homogeneous\/}, Bull. Amer. Math. Soc. 75 (1969), 554-558.

\bibitem{Fearnleyunique2} L. Fearnley, {\em The pseudo-circle is unique\/},
    Trans. Amer. Math. Soc. 149 (1970), 45-64.
    
\bibitem{Hagopian} C.L. Hagopian, {\em The fixed-point property for almost chainable homogeneous continua\/},
Illinois J. Math. 20 (1976),  650-652. 

\bibitem{Hu} S.T. Hu, Homotopy Theory, Elsevier Science and Technology Books, 1959.

\bibitem{KennedyRogers} J. Kennedy, and J.T. Rogers,  Jr., {\em Orbits of the pseudocircle\/},
Trans. Amer. Math. Soc. 296 (1986),  327-340. 

\bibitem{Knaster} B. Knaster, {\em Un continu dont tout so{\'u}s-continu est ind{\'e}composable\/}, Fund. Math. 3 (1922), 247-286.

\bibitem{Kuratowski} K. Kuratowski, Topology, Vol. II, Academic Press,  1968.

\bibitem{Lewis3} W. Lewis, {\em  Almost chainable homogeneous continua are chainable\/},
Houston J. Math. 7 (1981),  373-377. 

\bibitem{Lewis2} W. Lewis, {\em  The pseudo-arc\/}, Contemp. Math. 117 (1991), 103-123.

\bibitem{Lewis1} W. Lewis, {\em The pseudo-arc\/}, Bol. Soc. Mat. Mexicana 5 (1999), 25-77.

\bibitem{Moise} E.E.  Moise, {\em An indecomposable plane continuum which is homeomorphic to each of its nondegenerate subcontinua\/}, Trans. Amer. Math. Soc.  63 (1948), 581-594.

\bibitem{Rogers1} J.T. Rogers, Jr., {\em The pseudo-circle is not homogeneous\/}, Trans. Amer. Math. Soc. 148 (1970), 417-428.

\bibitem{Rogers2} J.T. Rogers, Jr., {\em Homogeneous,
separating plane continua are decomposable\/}, Michigan Math. J. 28 (1981), 317-322.

 \end{thebibliography}
\end{document}